\documentclass[12pt]{article}
\usepackage[dvipdfmx]{graphicx,color}
\usepackage{amsmath,amsfonts,amssymb}
\usepackage{hyperref}
\usepackage{ulem}

\usepackage[numbers]{natbib}
\newcommand{\sr}[1]{{\cal #1}}
\newcommand{\dd}[1]{\mathbb{#1}}

\newcommand{\eq}[1]{(\ref{eq:#1})}
\newcommand{\lem}[1]{Lemma~\ref{lem:#1}}

\newcommand{\thr}[1]{Theorem~\ref{thr:#1}}

\newcommand{\app}[1]{Appendix~\ref{app:#1}}
\newcommand{\sectn}[1]{Section~\ref{sect:#1}}

\newcommand{\thrt}[1]{\ref{thr:#1}}

\newcommand{\sect}[1]{\ref{sect:#1}}

\newcommand{\ol}{\overline}

\newcommand{\pend}{\hfill \thicklines \framebox(6.6,6.6)[l]{}}

\newenvironment{proof}{\noindent {\sc  Proof.} \rm}{\pend}
\newenvironment{proof*}[1]{\noindent {\sc  #1} \rm}{\pend}

\renewcommand{\theequation}{\thesection.\arabic{equation}}

\newtheorem{theorem}{Theorem}[section]
\newtheorem{lemma}{Lemma}[section]

\newtheorem{proposition}{Proposition}[section]

\newtheorem{remark}{Remark}[section]

\newcommand{\setsection}[2] {
\setcounter{section}{#1}
\setcounter{subsection}{0}
\setcounter{equation}{0}
\setcounter{conjecture}{0}
\setcounter{assumption}{0}
\setcounter{question}{0}
\setcounter{definition}{0}
\setcounter{theorem}{0}
\setcounter{corollary}{0}
\setcounter{lemma}{0}
\setcounter{proposition}{0}
\setcounter{remark}{0}
\setcounter{appen}{0}
\setsection*{\large \bf \thesection. #2}}

\newcommand{\setnewcounter} {
\setcounter{subsection}{0}
\setcounter{equation}{0}
\setcounter{conjecture}{0}
\setcounter{assumption}{0}
\setcounter{question}{0}
\setcounter{definition}{0}
\setcounter{theorem}{0}
\setcounter{corollary}{0}
\setcounter{lemma}{0}
\setcounter{proposition}{0}
\setcounter{remark}{0}
}

\begin{document}
\title{\bf \Large Tails in a fixed-point problem for a branching process with state-independent
immigration}

%%%%%%%%%%%%%%
%   AUTHORS  %
%%%%%%%%%%%%%%
\author{Sergey Foss\\Heriot-Watt University and \\ Novosibirsk State University\footnote{Research of S.~Foss is supported
by RSF research grant No. 17-11-01173}\\ \and Masakiyo Miyazawa\\ Tokyo University of Science\footnote{Research of M.~Miyazawa is supported by
JSPS KAKENHI Grant No. JP16H02786}}
%\date{\today, version $20180828-2$}
%\date{}

\maketitle

\begin{center}
{\bf The authors congratulate Guy Fayolle on the occasion of his\\ birthday
 and 
wish him many more happy and productive years!}
\end{center}

\begin{abstract}
We consider a fixed-point equation for a non-negative integer-valued random variable, that appears in branching processes with state-independent  immigration. A similar equation appears in the analysis of a single-server queue with a homogeneous Poisson input, feedback and permanent customer(s). 

It is known that the solution to this equation uniquely exists under mild first and logarithmic moments conditions. We find further the tail asymptotics of the distribution of the solution when the immigration size and branch size distributions are heavy-tailed. We assume that the distributions of interest are dominantly varying  and have a long tail. This class includes, in particular, (intermediate, extended) regularly varying distributions.

We consider also a number of generalisations of the model.
\end{abstract}
  
\begin{quotation}
\end{quotation}

\section{Introduction}
\label{sect:introduction}

We are interested in the following fixed-point equation for a non-negative integer-valued random variable $X$,
\begin{align}
\label{eq:fixed point 1}
  X =_{st} A + \sum_{i=1}^{X} B_{i},
\end{align}
where ``$=_{st}$'' represents the equality in distribution, $A$ and $B_{i}$ are nonnegative integer-valued random variables independent of $X$, and $B_{i}$ for $i=1,2,\ldots$ are i.i.d
random variables that do not depend on $A$. In what follows, we drop the subscript $i$ from $B_{i}$ if its distribution is only concerned.

Equation \eq{fixed point 1} arises in a branching process with state independent immigration, and its solution is subject to the stationary distribution of this branching process if it exists. Thus, the distribution of $X$ is the stationary distribution of a relatively simple discrete-time Markov chain with state space $\dd{Z}_{+} \equiv \{0,1,\ldots\}$. Although its existence and uniqueness are well studied (e.g., see \cite{FostWill1971}) and its moments can be inductively computed if they are finite, little is known about the distribution itself.  Its tail asymptotics have been studied recently in 
\cite{BasrKuliPalm2013} in the case of regularly varying distributions with a limited range of their parameters,  
and \cite{BarcBoszPap2018} have extended their results onto the case of two-level processes.

This is contrasted with a reflecting random walk on $\dd{Z}_{+}$, whose stationary distribution is the distributional solution of $X$ satisfying the following fixed point equation:
\begin{align}
\label{eq:fixed point 2}
  X =_{st} \max(0,X+\xi), 
\end{align}
where $X$ and $\xi$ are mutually independent. 
Various aspects of the solution of this fixed point equation and of related one- and multidimensional  problems have been studied in the queueing literature and, in particular, using the modern theory of   
random walks where Guy Fayolle and his colleagues have made a great 
contribution, see their books \cite{FayoIasnMaly1999,FayoMalyMens1995}. 

The fixed-point equation \eq{fixed point 1} is highly nonlinear compared with \eq{fixed point 2}, which causes a difficulty to see how does the distribution of $X$ look like. In this paper, we are focusing on the study of its tail asymptotics, with assuming $A$ and $B$ to have heavy-tailed distributions, that are both dominantly varying and long-tailed. For this, we introduce a new approach which allows to go beyond
the regular variation. In particular, in the regularly varying case, we extend the results of \cite{BasrKuliPalm2013}, with providing simpler proofs (see
our Theorems \thr{X 1} and \thr{X 2}, case (I)), and consider another situation (\thr{X 2}, case (II)).
% where a different type of the tail asymptotics does  appear. 
This paper is close
to \cite{FossMiya2018} where the tail asymptotics for related objects have been studied. 

There is another type of fixed point equations, which arise in internet page ranking problems (see \cite{AlsmMein2013,JeleOlve2012}) and in a queueing problem (see \cite{AsmuFoss2018}). The heavy tail asymptotics are also studied for their solutions. However, it should be noted that those fixed point equations are essentially different from \eq{fixed point 1}.

In what follows, we say that two (strictly) positive functions $f(x)$ and $g(x)$ are {\it asymptotically equivalent}
(at infinity) and write $f(x)\sim g(x)$ if $\lim_{x\to\infty} f(x)/g(x) =1$.
We write $f(x) \gtrsim g(x)$ if $\liminf_{x\to\infty} f(x)/g(x)\ge 1.$ For two random variables $X$ and $Y$, equality $X=_{st}Y$ means that $X$ and $Y$ are identically distributed.

We list below some known classes of heavy-tailed distributions which are used for a reference distribution $G$ either on 
$\dd{R}_{+}$ or on $\dd{R}$. We use the same notation $G$ for the distribution and for its distribution function. Let $\overline{G}(x)=1-G(x)$ be the tail distribution function. 
\begin{enumerate}
\item $G$ belongs to the class $\sr{L}$ of {\it long-tailed} distributions if, 
for some (or, equivalently, for any) $y$ and 
as $x\to\infty$, $\frac{\ol{G}(x+y)}{\ol{G}(x)} \rightarrow 1$ (we may write equivalently $\ol{G}(x+y) \sim \ol{G}(x)$).

\item $G$ belongs to the class $\sr{S}$ of {\it subexponential} distributions if $G\in \sr{L}$ and if $\ol{G*G}(x) \sim 2\ol{G}(x)$. For distributions on the positive half-line, the 2nd condition implies the 1st. 

\item $G$ belongs to the class $\sr{S}^{*}$ of {\it strong subexponential} distributions if $G$ belongs to $\sr{L}$, if 
$m_{G}^+ \equiv \int_{0}^{\infty} \ol{G}(x) dx$ is finite and if 
$\int_0^x \ol{G}(y) \ol{G}(x-y) dy \sim
2m_{G}^+ \ol{G}(x)$.% \quad \mbox{as}$. % for $x\to\infty$.

\item $G$ belongs to the class $\sr{D}$ of {\it dominantly varying} distributions  if there exists $y >1$ such that $\liminf_{x\to\infty}
\frac{\ol{G}(y x)}{\ol{G}(x)} > 0$. Then the same holds for any $y>1$. 

\item $G$ belongs to the class $\sr{IRV}$ of {\it intermediate regularly varying} distributions if
\begin{align}
\label{eq:IRV}
  \lim_{y \downarrow 1} \liminf_{x\to\infty}
\frac{\ol{G}(y x)}{
\ol{G}(x)} = 1.
\end{align}

\item $G$ belongs to the class $\sr{ERV}$ of {\it extended regularly varying} distributions if there are some $\alpha_{+}, \alpha_{-} > 0$ such that
\begin{align}
\label{eq:ERV}
  y^{-\alpha_{-}} \le \liminf_{x\to\infty} \frac{\ol{G}(y x)}{\ol{G}(x)} \le \limsup_{x\to\infty} \frac{\ol{G}(y x)}{\ol{G}(x)} \le y^{-\alpha_{+}}, \qquad \forall y \ge 1.
\end{align}

\item $G$ belongs to the class $\sr{RV}$ of {\it regularly varying} (at infinity) distributions if, for some 
$\alpha >0$, 
\begin{align}
\label{eq:RV}
  \ol{G}(x) = x^{-\alpha} L(x),
\end{align}
where $L(x)$ is a {\it slowly varying} (at infinity) function, i.e. $L(cx)\sim L(x)$, for any $c>1$. 
\end{enumerate}

A simple example of the tail distribution function from the class $\sr{ERV}\setminus\sr{RV}$ is given in \app{ERV}.
Note that the definitions of $\sr{ERV}$ and $\sr{RV}$ also can be used for positive valued functions on $[0,\infty]$ instead of $\ol{G}$. We recall some basic properties of heavy-tailed distributions we refer to in the paper. More properties and details may be found, e.g., in the books \cite{BingGoldTeug1987,FossKorsZach2011}. In particular, $\sr{ERV}$ is well studied in \cite{Clin1994,DrasSene1976}.

First, note the following relations between the introduced classes of heavy-tailed distributions: 
\begin{align}
\label{eq:class order}
  \sr{RV} \subset \sr{ERV} \subset \sr{IRV} \subset \sr{L} \cap \sr{D} \subset \sr{S}^* 
  \subset \sr{S} \subset \sr{L},
\end{align}
where each of the inclusions is strict.
Further, each of these classes in {\it closed} with respect to the tail equivalence: if $G_1$ belongs to a class,
and if $\overline{G}_1(x) \sim \overline{G}_2(x)$, then $G_2$ belongs to the same class. 
We also need the following property: if $G_1$ belongs to one of the classes $\sr{RV}, \sr{IRV}, \sr{L} \cap \sr{D}$, $\sr{S}^*$ or $\sr{S}$  
and if $\overline{G_2}(x) = o(\overline{G_1}(x))$ as $x\to\infty$,
then, for any fixed $j=1,2,\ldots$, we have $\overline{G_1*G_2^{(*j)}}(x) \sim \overline{G_1}(x)$ (we may informally say that the tail of $G_2^{(j)}$ is ``ignorable'') and, therefore, $G_1*G_2^{(*j)}$ belongs to the same class.   

The following result will be repeatedly used in our proofs.
\begin{proposition}\label{prop1} 
{\bf (Theorem 1 in \cite{FossZach2003})}\\
Let $S_n=\sum_1^n \xi_i$, $n=1,2,\ldots$ be the sums of i.i.d.r.v.'s with negative mean, $S_0=0$ and 
$M_n = \max_{0\le k \le n} S_n$. Let $\sigma\le \infty$ be any stopping time with respect to $\{\xi_n\}$.
If the common distribution $G$ of the $\xi$'s is strong subexponential, then
\begin{align*}
\lim_{x\to\infty} \frac{\dd{P}(M_{\sigma}>x)}{\overline{G}(x)} = \dd{E} \sigma \le \infty.
\end{align*}
\end{proposition}

\section{The tail asymptotics for the solution of the fixed point equation}\setnewcounter
\label{sect:solution}

As we have already said, it is hard to derive an analytically tractable solution of the distribution of $X$, so we consider its tail asymptotics. To exclude a trivial exception, we assume throughout the paper that
\begin{align}
\label{eq:A 0}
  0 < \dd{P}(A = 0) < 1, \qquad \dd{P}(B = 0) < 1.
\end{align}
Let $a = \dd{E}(A)$ and $b = \dd{E}(B)$. The following fact is well-known.

\begin{lemma}
\label{lem:stability 1}
{\bf (Theorem and Corollary 2 of \cite{FostWill1971} and Theorem 3.1 of \cite{Sene1971})}\\ 
$(I)$ If $b > 1$, then $\eq{fixed point 1}$ has no solution.\\
$(II)$ If $0 < b \le 1$, then the solution $X$ of $\eq{fixed point 1}$ has a proper distribution if and only if
\begin{align}
\label{eq:stability 1}
  \int_{0}^{1} \frac {1-\dd{E}(x^{A})} {\dd{E}(x^{B}) - x} dx < \infty.
\end{align}
In particular, if $0 < b < 1$, then the condition $\eq{stability 1}$ can be replaced by
\begin{align}
\label{eq:stability 2}
  \dd{E} (\log \max(A,1)) < \infty,
\end{align}
which obviously holds if $a < \infty$, and the solution $X$ is unique in distribution.
\end{lemma}

\begin{remark}
\label{rem:stability 1}
Seneta (\cite{Sene1971}, Theorem 3.1) proved uniqueness of the distribution of $X$ in terms of generating functions in Theorem 3.1. In particular, he showed that the probability generating function $\dd{E}(s^{X})$ is regularly varying as $s \uparrow 1$. % when $\dd{P}(B < \infty) = 1$, which we always assume.
\end{remark}

In this paper, we assume that the stability conditions $b<1$ and \eq{stability 2}, and there is a reference distribution $G$ such that
\begin{align}\label{eq:c1c2}
  \lim_{x \to \infty} \frac {\dd{P}(A>x)} {\ol{G}(x)} = c_{1}, \qquad \lim_{x \to \infty} \frac {\dd{P}(B>x)} {\ol{G}(x)} = c_{2},
\end{align}
for some constants $c_{1},c_{2} \ge 0$ such that $c_{1}+c_{2} >0$. We consider the following three cases: 
\begin{itemize}
\item [(i)] $c_{1}, c_{2} > 0$. Here $\dd{P}(A>x) \sim \frac {c_{1}}{c_{2}} \dd{P}(B>x)$.
\item [(ii)] $c_{1} = 0, c_{2} > 0$. Here $\dd{P}(A>x) = o(\dd{P}(B>x))$.
\item [(iii)] $c_{1} > 0, c_{2} = 0$. Here $\dd{P}(B>x) = o(\dd{P}(A>x))$.
\end{itemize}
Note that if $b$ is finite, then $a$ must be finite in the cases (i) and (ii). However, in the case (iii), $a$ is allowed to be either finite or infinite. We consider these situations separately.

\begin{remark}
\label{rem:00}
In what follows, we consider the ``sub-critical'' case only: $0<b<1$. 
In this case, the two tail asymptotics, of $\dd{P}(A>x)$ and of $\dd{P}(B_1>x)$,
may contribute to that of $\dd{P}(X>x)$. The case $b=1$ is more involved, and it is a subject of a future research. In this case, the tail of a null-recurrent random walk may also contribute to the asymptotics of $\dd{P}(X>x)$. 
\end{remark}

We need another condition on the distribution $G$. For $c>1$, let
$$
T_c(x) = \sum_{n=0}^{\infty} \ol{G}(c^n x), \qquad x > 0.
$$
We note the following fact.
\begin{lemma}
\label{lem:T_c}
If $b < 1$ and \eq{stability 2} holds, then $T_c(x)$ is finite for all $x > 0$ and $c > 1$.

\begin{proof}
$T_c(x)$ is obviously finite if $G$ has a finite mean. Otherwise, we must have (iii) because $b < 1$, then there is an $x_{0}>1$ for each $\epsilon > 0$ such that $c_{1} \ol{G}(x) \le (1+\varepsilon) \dd{P}(A > x) $ for each $x \ge x_{0}$, and therefore, by \eq{stability 2},
\begin{align*}
  T_{c}(x) \le \frac {1+\varepsilon} {c_{1}} \sum_{n=0}^{\infty} \dd{P}(A > c^{n} x) \le \frac {1+\varepsilon} {c_{1}} \sum_{n=0}^{\infty} \dd{P}(\log \max(A,1) > n \log c + \log x) < \infty.
\end{align*}
For each $x > 0$, we choose $n_{0}$ such that $c^{n_{0}}x \ge x_{0}$, then
\begin{align*}
  T_{c}(x) = \sum_{n=0}^{n_{0}-1} \ol{G}(c^n x) + \sum_{n=n_{0}}^{\infty} \ol{G}(c^{n-n_{0}}c^{n_{0}}x) \le \sum_{n=0}^{n_{0}-1} \ol{G}(c^n x) + \sum_{n=0}^{\infty} \ol{G}(c^{n} x_{0}) < \infty.
\end{align*}
 Hence,  the lemma is proved.
\end{proof}
\end{lemma}

 Then the extra condition is: with $c_0=1/b$, 
\begin{align}
\label{eq:G 22}
\lim_{c \uparrow c_0} \limsup_{x\to\infty} T_c(x)/T_{c_0}(x) =
\lim_{c\downarrow c_0} \liminf_{x\to\infty} T_c(x)/T_{c_0}(x) = 1.
\end{align}
The condition \eq{G 22} is not easy to check. So, we provide below a sufficient condition for it to hold.

\begin{lemma}
\label{lem:G22 sufficient}
The condition \eq{G 22} is satisfied if $G \in \sr{IRV}$ and if $\ol{G}$ has the Karamata upper index $c_{+}(\ol{G}) < 0$, that is, 
\begin{align*}
  c_{+}(\ol{G}) \equiv \inf \left\{c \in \dd{R}; \limsup_{x \to \infty} \frac{\ol{G}(\lambda x)} {\ol{G}(x)} \le \lambda^{c} \mbox{ uniformly in } \lambda \in [1,\Lambda], \mbox{ for all} 
  \ \Lambda > 0 \right\} < 0.
\end{align*}
In particular, if $G \in \sr{ERV}$, then \eq{G 22} holds, and the function $T_{1/b}(x)$ also belongs to the class $\sr{ERV}$ of non-increasing functions.
\end{lemma}

We prove this lemma in \app{G22 sufficient}. An extended regularly varying distribution is a special case of the distributions with $c_{+}(\ol{G}) < 0$ (see, e.g., \cite[Theorem 2.1.1]{BingGoldTeug1987}). However, we do not know in general how large is the subclass of distributions from ${\cal L\cap D}$ that satisfy \eq{G 22}, and
this is a subject of further studies.  

\begin{theorem}
\label{thr:X 1}
Assume that ${\eq{A 0}}$ holds and that $a<\infty$ and $0 < b \equiv \dd{E}(B) < 1$. Assume that either one of conditions (i), (ii) or (iii)  holds.

If $G$ is strong subexponential, then, for any $d_2>1/b$, 
\begin{align}
\label{eq: new1}
\dd{P}\left( X > x \right) \gtrsim \frac {(1-b) c_{1} + a c_{2}}{1-b}T_{d_2}(x).
\end{align}

Further, if  $G\in \sr{D} \cap \sr{L}$, then, for any $1<d_1<1/b$, 
\begin{align}
\label{eq: new2}
\dd{P}\left( X > x \right) \lesssim \frac {(1-b) c_{1} + a c_{2}}{1-b}T_{d_1}(x).
\end{align}
 
If, in addition, \eq{G 22} is satisfied, then
\begin{align}
\label{eq:X 2}
  \dd{P}\left( X > x \right) \sim \frac {(1-b) c_{1} + a c_{2}}{1-b} T_{1/b}(x).
  % \sum_{n=0}^{\infty} \ol{G}(b^{-n} x).%, \quad  \mbox{as} \ x\to\infty.
\end{align}
%Here we use the convention that $a c_{2} = 0$ for $c_{2} = 0$ and $a=\infty$.

In particular, \eq{X 2} holds if $G$ has an extended regularly varying distribution $\eq{ERV}$. In this case, $X$ also has an extended regularly varying distribution. Furthermore, if $G$ has a regularly varying distribution $\eq{RV}$, then $X$ also has a regularly varying distribution with the same parameter $\alpha$, and
\begin{align}
\label{eq:Xi-2}
  \dd{P}\left( X > x \right) \sim \frac {(1-b) c_{1} + a c_{2}} {(1-b)(1-b^{\alpha})} \ol{G}(x).%, \quad \mbox{as } \ x \to \infty.
\end{align}
\end{theorem}

\begin{theorem}
\label{thr:X 2}
Assume again conditions ${\eq{A 0}}$ and $0 < b \equiv \dd{E}(B) < 1$ to hold.
Assume that $G \in \sr{D} \cap \sr{L}$, and that condition (iii)  holds.
Assume now that $a=\infty$, but that \eq{stability 2} holds. 
Assume further that either \\
$(I)$ random variable $B$ has finite variance $\sigma^2 = {\rm Var} B$ and
$C:= \liminf_{x\to\infty}x\ol{G}(x)\in (0,\infty]$ or\\
$(II)$ 
the distribution 
$H_{I}(x) \equiv 1- \min(1,\int_{x}^{\infty} \dd{P}(B > u) du)$ 
is subexponential and   
\begin{align}
\label{eq:c2 additional}
\limsup_{x \to \infty} \ol{H}_{I}(x)/\ol{G}(x)  < \infty.
\end{align}
If \eq{G 22} is satisfied, then the asymptotic \eq{X 2} holds again (here we use the convention that $a c_{2} = 0$ for $c_{2} = 0$ and $a=\infty$).
In particular, this is the case for $G \in \sr{ERV}$, and therefore \eq{Xi-2} holds for $G \in \sr{RV}$.
\end{theorem}

We prove these theorems in \sectn{proof}.

\begin{remark}
\label{rem:Z1}
The asymptotics $\eq{X 2}$ hold if and only if  
%equivalent to 
%If \eq{G 1} is not satisfied in Theorem \ref{thr:X 1}, we only have
\begin{align}
\label{eq:Xi-1}
 \dd{P}\left( x < X \le x/b \right) \sim \frac {(1-b) c_{1} + a c_{2}}{1-b} \ol{G}(x).
 \end{align}
%If \eq{G 1} holds, then both sides of \eq{Xi-1} are asymptotically equal, %which is equivalent to \eq{X 2}. 
This follows from the fact that $T_{1/b}(x) - T_{1/b}(x/b) = \ol{G}(x)$. Let $D = \frac {(1-b) c_{1} + a c_{2}}{1-b}$, then $\eq{X 2}$ implies that
\begin{align*}
  \dd{P}(X>x) - \dd{P}(X>x/b) = D (T_{1/b}(x) - T_{1/b}(x/b)) (1+ o(1)) = D \ol{G}(x) (1+o(1)),
\end{align*}
which proves $\eq{Xi-1}$. On the other hand, if $\eq{Xi-1}$ holds, then summing up $\eq{Xi-1}$ by replacing $x$ by $x/b^{n}$ and applying $\eq{Xi-1}$ to each term yield $\eq{X 2}$. It is notable that the condition $\eq{G 22}$ is not needed for the equivalence of $\eq{X 2}$ and $\eq{Xi-1}$.
\end{remark}

\begin{remark}
\label{rem:Z2}
The asymptotics $\eq{X 2}$ have a natural interpretation in the terms of the Principle of a Single Big Jump (PSBJ): for the sum to be large, either one of the summands or the counting random variable must be large. Namely, we may rewrite equation $\eq{fixed point 1}$ as an a.s. equation:
\begin{align}
\label{eq:fixed point 111}
  X^{(1)} = A + \sum_{i=1}^{X^{(2)}} B_{i},
\end{align}
where r.v.'s $X^{(1)}$ and $X^{(2)}$ have the same distribution, and all r.v.'s on the right are mutually independent. 
Then $X^{(1)}>x$ if either $A>x$ or one of $B_i>x$, for $i=1,\ldots, X^{(2)}$, or if $X^{(2)}$ is ``sufficiently large''.
The latter means that, due to the Strong Law of
Large Numbers, we have to have $\sum_1^{X^{(2)}} B_i \approx bX^{(2)} >x$, which leads to $X^{(2)}>x/b$. Therefore,
\begin{align*}
\dd{P}(X^{(1)}>x) \approx \dd{P} (A>x) + \dd{E} X^{(2)} \cdot \dd{P}(B_1>x) + \dd{P}(X^{(2)}>x/b).
\end{align*}
Then the induction argument completes the derivation. However, this is not a proof, and just an explanation of the phenomenon. 
\end{remark}

\begin{remark}\label{rem:Z3}
Note that the two conditions $a=\infty$ and
$\liminf_{x\to\infty}x\ol{G}(x) >0$ do hold if $G$ is a regularly varying  distribution with parameter $\alpha \in [0,1]$.
So, the case (I) of our Theorem \ref{thr:X 2} extends the corresponding result
of \cite{BasrKuliPalm2013}. In addition, these two conditions may hold for
other distributions from the class ${\cal D}\cap {\cal L}$ which are not regularly
varying -- for example, for distributions from the class ${\cal ERV}$. 
\end{remark}

\begin{remark}\label{rem:open}
We do not know how restrictive are the conditions of \thr{X 2},
% and whether is \eq{X 2} the only asymptotics one can expect in the class 
%$\sr{L}\cap\sr{D}$ of heavy-tailed distributions. See 
see Section \sect{OP} for related comments.
\end{remark}

%\begin{remark}
%\label{rem:X 1}
%Since $\sr{IRV} \subset \sr{D} \cap \sr{L}$, \eq{X 2} holds if $G$ is intermediate regularly varying.
%\end{remark}

In the rest of this section, we discuss stochastic models where the fixed-point equation \eq{fixed point 1} arises. In what follows, we assume \eq{A 0}, \eq{stability 2} and that $0 < b < 1$, $a < \infty$, and therefore \eq{fixed point 1} has a unique distributional solution, by \lem{stability 1}.
 
Define $X_{n}$ inductively for $n=0,1,\ldots$ by
\begin{align}
\label{eq:Markov X 1}
  X_{n} = \left\{
\begin{array}{ll}
 0, & n = 0,\\
 A_{n} + \sum_{i=1}^{X_{n-1}} B_{i,n}, \qquad & n \ge 1,
\end{array}
\right.\end{align}
where $A_{n}$ and $B_{i,n}$ have the same distributions as $A$ and $B$, respectively, and $\{A_{n}\}$ and $\{B_{i,n}\}$ are sequences of $i.i.d.$ random variables. Clearly, $\{X_{n}; n=0,1,\ldots\}$ is a Markov chain with state space $\dd{Z}_{+}$, and a branching process with immigration $\{A_{n}\}$. The following lemma is a direct consequence of \eq{Markov X 1} and \lem{stability 1}.

\begin{lemma}
\label{lem:monotone 1}
Let the Markov chain $\{X_{n}; n=0,1,\ldots\}$ be defined by $\eq{Markov X 1}$.\\ 
$(I)$ The distribution of $X_{n}$ is stochastically non-decreasing in $n$, i.e. $\dd{P}(X_{n+1}>x) \ge \dd{P}(X_n>x)$, for all $n$ and $x$, and $X_{n}$ converges to the solution $X$ of $\eq{fixed point 1}$ in distribution as $n \to \infty$.\\
$(II)$ If $0 < b < 1$ and $a < \infty$, then $\dd{E}(X_{n})$ is non-decreasing in $n$ and
\begin{align}
\label{eq:EX 1}
  \lim_{n \to \infty} \dd{E}(X_{n}) = \frac a{1-b} < \infty,
\end{align}
and therefore $\dd{E}(X) = a/(1-b) < \infty$ for the solution $X$ of $\eq{fixed point 1}$, by the uniqueness of \lem{stability 1}.
\end{lemma}
For completeness, we provide a proof of this lemma in \app{G23}.

The branching process $\{X_{n};n=0,1,\ldots\}$ also arises in a feedback single server queue with Poisson arrivals and with one (or more) permanent customer(s). 
In this model, all the arriving customers %, except the ``target customer'', 
receive service in the First-Come-First-Served order.
There are two types of customers, the first (``target'') customer and all other customers. 
The target customer
arrives at time $0$ at the empty system and, after each service completion, with probability one joins the end of the queue again, for another service. Any other customer, after its service completion, either joins again the end of the queue (for another service), with probability $p$, or leaves the system, with probability $1-p$. 
If we denote by $X_{n}$ the number of customers observed by the tagged customer at its $n$-th return to the queue, then
a sequence $\{X_n\}$ forms a time-homogeneous Markov chain which converges to the stationary distribution (under natural
assumptions), and a random variable $X$ with that distribution satisfies the fixed-point equation. See Section 2.2 of \cite{FossMiya2018} for related details. 
%Some related model is considered in \sectn{multiple}. This $X_{n}$ plays a crucial role to compute the sojourn time of a feedback single server queue in \cite{FossMiya2018}. We may replace the tagged customer by a permanent customer who stays in the system forever, then $X_{n}$ is well defined for all $n \ge 0$, and agrees with that of the tagged customer returning to the queue $n$ times. Thus, \thr{X 1} is also interesting in the queueing model context, although the permanent customer is not considered in \cite{FossMiya2018}.

%\begin{example}
%\label{exa:multiple permanent customers}
Similarly, we may consider a sub-critical branching process with ${k}$ permanent particles. We may rephrase a particle as a customer. Each ordinary customer produces (independently of
everything else) a random number of offspring 
with distribution $G$ and either stays in the system, with probability $p$, or leaves the system (dies),
with probability $q=1-p$,  while 
each permanent customer produces (again independently of everything else) a random number of offspring with distribution $G$ and stays
in the system. (Again, one can view this system as a single-server queueing system with Poisson input stream, where customers are served one-by-one and where customers arriving during service of
any customer are viewed as its ``offspring'').

Let $\xi$ be a random variable with distribution $G$, and let $Y$ be the population in the steady state.
Here we assume that $\dd{E}(\xi) + p<1$ (to have subcriticality),  and the fixed-point equation
looks like
\begin{align}
\label{eq:multiple 1}
  Y=_{st} {k} + \sum_{i=1}^{Y-{k}} \alpha_i + \sum_{i=1}^Y \xi_i,
\end{align}
where $\alpha_i$ are Bernoulli random variables with parameter $p$, random variables $\xi_{i}$ have distribution $G$, and all random variables on the right of equation \eq{multiple 1} are mutually independent. Let $X = Y - k$, and let
\begin{align*}
  A = \sum_{i=1}^{k} \xi_{i,1}, \qquad B_{i} = \alpha_i+\xi_{i,2},
\end{align*}
where $\xi_{i,1}$ and $\xi_{i,2}$ are $i.i.d.$ and distributed as $\xi$, then \eq{multiple 1} is the exactly same form of a fixed point equation as \eq{fixed point 1}. Thus, we can get the tail asymptotic \eq{X 2} for $X$ and therefore for $Y = X + k$ under the assumptions of \thr{X 1}.
%\end{example}

\section{Proof of Theorems \thrt{X 1} and \thrt{X 2}}\setnewcounter
\label{sect:proof}

{\sc Proof of \thr{X 1}}.
We first consider the case (i). Here since $b$ is finite, $a$ is finite too, because of the tail equivalence. Hence, $\dd{E}(X) <\infty$ by \lem{monotone 1}. We also note that $G$ is strongly subexponential by \eq{class order}.
For $0 < y < x$, let
\begin{align*}
  & I^{-}(x,y) = \dd{P}\left( A + \sum_{i=1}^{X} B_{i} > x, X \le y \right),\\
  & I^{+}(x,y) = \dd{P}\left( A + \sum_{i=1}^{X} B_{i} > x, X > y \right),
\end{align*}
then the fixed-point equation \eq{fixed point 1} may be written as
\begin{align*}
  \dd{P}(X > x) & = I^{-}(x,y) + I^{+}(x,y).
\end{align*}
%We consider upper and lower bounds for $I^{-}(x,y)$ and $I^{+}(x,y)$.

First, we find the tail asymptotics for $I^{-}(x,y)$ as $x\to\infty$, for any fixed $y$. %Since a subexponential distribution is 
Here we need only the strong subexponentiality assumption.
We
may apply Theorem 3.37 of \cite{FossKorsZach2011} and obtain that, as $x\to\infty$,
\begin{align}
\label{eq:I- 1}
I^-(x,y) & = \dd{P}\left(A{\bf 1}{\{X\le y\}} + 
\sum_{i=1}^{X {\bf 1}{\{X\le y\}}} B_{i} > x\right)\nonumber \\
& = \sum_{j=0}^y \dd{P} (X=j) \dd{P} \left(
A + \sum_{i=1}^j B_i>x\right) \nonumber\\
& \sim \sum_{j=0}^y \dd{P}(X=j) (c_1+c_2j) \ol{G}(x) \nonumber\\
& \sim \dd{E}((c_{1} + c_{2} X){\bf 1}{\{X \le y\}}) \ol{G}(x).
\end{align}
Here ${\bf 1}(E)$ is the indicator of event $E$, it takes value $1$ if the event occurs and $0$, otherwise. 

We next establish the lower and upper bounds for $I^{+}(x,y)$. 
We start with the upper bound. Here we assume in addition that $G\in \sr{D} \cap \sr{L}$.
 We choose a sufficiently small $\varepsilon > 0$ such that $0 < b + \varepsilon < 1$. 
Consider a random walk $S_{n} = \sum_{i=1}^{n} (B_{i} - (b+\varepsilon/2))$ with initial value $S_0=0$. 
It has a negative drift: 
\begin{align*}
   \dd{E}B_{i} - (b+\varepsilon/2) = - \varepsilon/2 < 0.
\end{align*}
Let $d_1 = (b + \varepsilon)^{-1}$, then 
\begin{align}
\label{eq:I+ 1}
  I^{+}(x,y) & \le \dd{P}\left( X > d_{1} x \right) + \dd{P}\left( A + \sum_{i=1}^{X} B_{i} > x, y < X \le d_1 x \right).
\end{align}
In turn, the second term in the RHS may be bounded above by 
\begin{align}
\label{eq:I+ 2}
 & \quad \dd{P}\left( A + \sum_{i=1}^{X} B_{i} > x, y < X \le d_{1} x \right) \nonumber\\
 & \quad \le \dd{P}\left( A + S_{X} + (b+\varepsilon/2) X > x, y < X \le d_{1} x \right) \nonumber\\
 & \quad \le \dd{P}\left( A + \max_{1 \le n \le X} S_{n} + (b+\varepsilon/2) d_{1}x > x, y < X \le d_{1}x \right) \nonumber\\
 & \quad \le \dd{P}\left(A + \max_{1 \le n \le X} S_{n} > \frac {\varepsilon/2}{b+\varepsilon} x, X > y \right) \nonumber\\
 & \quad \le \dd{P}\left( A{\bf 1}\{X > y\} + \max_{1 \le n \le X {\bf 1}
 \{X > y\}} S_{n} > \frac {\varepsilon/2}{b+\varepsilon} x \right) \nonumber\\
 & \quad = \dd{P}(X>y) \dd{P} \left( A+ \max_{1\le n \le \widetilde{X}} S_n > \frac {\varepsilon/2}{b+\varepsilon} x \right) \nonumber\\
 & \quad \le \dd{P}(X>y) \left(\dd{P} \left( A > \frac {\varepsilon/4}{b+\varepsilon} x\right) + \dd{P} \left( \max_{1\le n \le \widetilde{X}} S_n > \frac {\varepsilon/4}{b+\varepsilon} x \right) \right)\nonumber\\
 & \quad \sim \dd{P}(X>y) \dd{E}\left( c_{1} + c_{2} \widetilde{X}\right) \ol{G}\left(\frac {\varepsilon/4}{b+\varepsilon} x \right)\nonumber\\
 & \quad = \dd{E}\left( (c_{1} + c_{2} X) {\bf 1}\{X > y\} \right) \ol{G}\left(\frac {\varepsilon/4}{b+\varepsilon} x \right). 
\end{align}
Here $\widetilde{X}$ has the conditional distribution $\dd{P}(\widetilde{X}>t) = \dd{P}(X>t \ | \ X>y)$ and does not depend on $A$ and $B_i$; and the equivalence in the pre-last line follows from Proposition \ref{prop1} since $G$ is strong subexponential. 

Since $G$ belongs to $\sr{D}$, there exist $\delta_{0}(\varepsilon)$ and $x_{0}(\varepsilon)$ such that
\begin{align*}
  \ol{G}\left(\frac {\varepsilon/4}{b+\varepsilon} x \right) \le \delta_{0}(\varepsilon) \ol{G}(x), \qquad \forall x \ge x_{0}(\varepsilon).
\end{align*}
Since $\dd{E}(X) < \infty$, one can choose sufficiently large $y=y_{\varepsilon,\delta_{1}}$ for any small $\delta_{1} > 0$ such that
\begin{align*}
  \dd{E}\left( (c_{1} + c_{2} X) {\bf 1}\{X > y\} \right) \delta_{0}(\varepsilon) \le \delta_{1},
\end{align*}
and therefore
so large that, for all $x \ge x_{0}(\varepsilon)$ and $y \ge y_{\varepsilon,\delta_{1}}$,
\begin{align*}
\dd{E}\left( (c_{1} + c_{2} X) {\bf 1}\{X > y\} \right) \ol{G}\left(\frac {\varepsilon/4}{b+\varepsilon} x \right)
&
\le \delta_{1} \ol{G}(x).
\end{align*}
Combining this inequality with \eq{I- 1} and \eq{I+ 1}, we conclude that there exists $x_0 \equiv x_0(\varepsilon,\delta_{1})$ such that, for any $x\ge x_0$,
\begin{align}
\label{eq:upper X 1}
  \dd{P}( x < X \le d_{1} x ) & \sim \dd{E}((c_{1} + c_{2} X){\bf 1}{\{X \le y\}}) \ol{G}(x) + \dd{P}\left( A + \sum_{i=1}^{X} B_{i} > x, y < X \le d_{1} x \right) \nonumber\\
 & \lesssim (\delta_{1} + c_{1}+ c_{2} \dd{E}(X)) \ol{G}(x).
\end{align}
%By the condition \eq{G 1}, this implies that, for some constant $C > 0$ 
%and $\forall \varepsilon \in (0, \varepsilon_{0}]$ for the $\varepsilon_{0}$ %of \eq{G 1} and for any sufficiently large $x$,
This implies that, for all sufficiently large $x$,
\begin{align}
\label{eq:upper X 2}
\dd{P}(X>x) & = \sum_{k=0}^{\infty} \dd{P} (xd_{1}^k < X \le xd_1^{k+1}) \nonumber\\
& \le (\delta_{1} + c_{1}+ c_{2} \dd{E}(X)) T_{d_1}(x)
%\sum_{k=0}^{\infty} C \ol{G}(d_{1}^{k-1} x) \le
%(\delta_{1} + c_{1}+ c_{2} \dd{E}(X)) C \ol{G}(x).
\end{align}  

Therefore, we get:
\begin{align*}
\limsup_{x\to\infty}\dd{P}(X>x)/T_{d_1}(x) & \le (\delta_{1} + c_{1}+ c_{2} \dd{E}(X)) ,
\end{align*}
for any sufficiently small positive $\varepsilon$ and $\delta_{1}$. Letting 
$\delta$ tend to zero, we obtain the following result:
\begin{align}
\label{eq:upper reg 1}
\limsup_{x\to\infty}\dd{P}(X>x)/T_{d_1}(x) & \le (c_{1} + c_{2} \dd{E}(X)).
\end{align}
 In particular, if condition \eq{G 22} holds, then we may tend $d_2$ to $1/b$
 and obtain the desired upper bound. Further, 
 if $G$ has an extended regularly varying distribution, then we have this upper bound by \lem{G22 sufficient}. This implies that, if $G$ has a regularly varying distribution with index $\alpha > 0$, 
 then we get the upper bound in \eq{Xi-2}. 
%\begin{align*}
%  C_{0,\alpha} = \sum_{k\ge 0} b^{k\alpha}.
%\end{align*}

We next consider the lower bound for $I^{+}(x,y)$. Letting $d_{2} = (b - \varepsilon)^{-1} > 0$ for a sufficiently small $\varepsilon > 0$, we have, for $y < x$,
\begin{align}
\label{eq:I+ 3}
  I^{+}(x,y) & \ge \dd{P}\left( A + \sum_{i=1}^{X} B_{i} > x, X > d_{2} x \right) \nonumber\\
  & \ge \dd{P}\left( X > d_{2} x \right) - \dd{P}\left( \sum_{i=1}^{X} B_{i} \le x, X > d_{2} x \right) \nonumber\\
  & \ge \dd{P}\left( X > d_{2} x \right) - \dd{P}\left( \sum_{i=1}^{d_{2} x} B_{i} \le x, X > d_{2} x \right) \nonumber\\
  & = \dd{P}\left( X > d_{2} x \right) \left(1 - \dd{P}\left( \sum_{i=1}^{d_{2} x} B_{i} \le (b-\varepsilon) d_{2} x\right)\right).
\end{align}
Here the subtrahend in \eq{I+ 3} decays exponentially fast, due to the
Chernoff's bound: since $B_i$ are positive, 
there exist universal positive constants $K$ and $\alpha$ such that, for all
$x>0$, 
\begin{align*}
  %\lim_{x \to \infty} 
  \dd{P}\left( \sum_{i=1}^{d_{2} x} B_{i} \le (b-\varepsilon) d_{2} x\right)
  \le
  Ke^{-\alpha x}. 
\end{align*}

Therefore, combining \eq{I- 1} with \eq{I+ 3} for any small $\varepsilon >0$ and $\delta_{2} >0$ , one can choose $y=y_{\varepsilon,\delta_{2}}$ so large that   
\begin{align}\label{eq:lower X 1}
  \dd{P}( x < X \le d_{2} x ) \ge (1-\delta_{2})(c_{1} + c_{2} \dd{E}(X)) \ol{G}(x) - Ke^{-\alpha x}, 
\end{align}
for all sufficiently large $x$.
Then, for the appropriate $K_1>K$, 
\begin{align*}
\dd{P}(X>x) & = \sum_{k=0}^{\infty} \dd{P} (xd_{2}^k < X \le xd_{2}^{k+1}) \ge
(1-\delta_{2})(c_{1} + c_{2} \dd{E}(X))T_{d_2}(x) - K_1e^{-\alpha x}\\
&\gtrsim (1-\delta_{2})(c_{1} + c_{2} \dd{E}(X))T_{d_2}(x), 
\end{align*}
for any small $\delta_2>0$. Letting $\delta_2$ tend to 0, we obtain the
desired lower bound.
Then, under assumption \eq{G 22}, we may let $d_2$ tend to $1/b$ and obtain
the lower bound that coincides with the upper bound obtained earlier. 
By \lem{G22 sufficient}, the statements in the last paragraph of this theorem are legitimated.
This completes the proof of the theorem in the case (i). 

The proof in the case (ii) is similar to the proof above, and even simpler. Since $\dd{P}(X\ge y)>0$, for any positive $y$, the distributional tail of $A$ is negligible with respect to that of $\sum_{i=1}^{X{\bf 1}(X\le y)} B_{i}$, of $\sum_{i=1}^{X{\bf 1}(X > y)} B_{i}$ (see Section 1 for the corresponding property), 
 and of $\max \{ S_{n}, 1 \le n \le X {\bf 1}(X > y)\}$. Then, clearly, \eq{I- 1} holds with $c_1=0$. Further, \eq{I+ 2} is also valid with $c_1=0$.  Then we get \eq{Xi-1} and \eq{Xi-2} for $c_{1} = 0$, and the proof is complete.

We next consider the case (iii). 
First, for any $j=1,2,\ldots$,
the tail distribution of $\sum_{i=1}^{j} B_{i}$ is negligible with respect to that of $A$ (see again Section 1). Therefore, 
we may take $y$ such that $\dd{P}(X\le y)>0$ and get:
\begin{align*}
I^-(x,y) & = \dd{P}\left(A {\bf 1}(X\le y) + \sum_{i=1}^{X {\bf 1}(X\le y)} B_{i} > x\right)\\
& = \sum_{j=0}^y \dd{P}(X=j) \dd{P} \left(A + \sum_{i=1}^j B_i>x\right)\\
& \sim \sum_{j=0}^y \dd{P}(X=j) \dd{P} (A >x)\\
& \sim c_{1} \dd{P}(X \le y) \ol{G}(x).
\end{align*}
For $I^{+}(x,y)$, we  use the same arguments as in the case (i), and then \eq{X 2} follows, with $a c_{2} = 0$. 

This completes the proof of \thr{X 1}.

{\sc Proof of \thr{X 2}}.
%Proof of (I). 
We again obtain the upper and lower bounds for $I^{-}(x,y)$ like in the case (i). 
However, since $a = \infty$, $\dd{E}(X)$ is infinite too, and we can not use the last two formulas of \eq{I+ 2} for getting the upper bound for $I^{+}(x,y)$ because $\dd{E}(X)$ is infinite. Therefore we modify these lines as follows.\\
In the case (I), we get
\begin{align}
\label{eq:I+ 2aa}
 & \quad \dd{P}\left( A + \sum_{i=1}^{X} B_{i} > x, y < X \le d_{1} x \right) \le 
 \dd{P}\left( A + \sum_{i=1}^{d_1x}B_i >x, y< X\le d_1x\right) \nonumber\\
 & \quad \le \dd{P}\left( A + \sum_{i=1}^{d_1x}B_i >x, y< X\right) = \dd{P}\left(A +  \sum_{i=1}^{d_1x}(B_i-b)> x(1-d_1b)\right)\dd{P}(X > y) \nonumber\\
 & \quad \le \left(\dd{P}(A>(1-d_1b)x/2) +
 \dd{P} \left( \sum_{i=1}^{d_1x} (B_i-b) >x (1-d_1b)/2 \right)\right) \dd{P}(X > y) \nonumber\\
 & \quad \le \left((c_1+o(1)) \ol{G}\left(\frac {\varepsilon}{2(b+\varepsilon)} x\right) + \frac {4(b+\varepsilon)^{2} d_{1} \sigma^{2}}{x\varepsilon^{2}}  \right) \dd{P}(X>y) 
\end{align}
where the inequality in the last line follows from the Chebyshev's inequality. Since $G \in \sr{D}$ implies that there is a $\delta(\varepsilon) > 0$ for each $\varepsilon > 0$ such that $\ol{G}\left(\frac {\varepsilon}{2(b+\varepsilon)} x\right) \le (\delta(\epsilon) + o(1)) \ol{G}(x)$ and since the condition (I) implies that $(C+ o(1)) x^{-1} \le \ol{G}(x)$, the last line of \eq{I+ 2aa} is not less than
\begin{align*}
 \left(c_1 \delta(\varepsilon) + o(1) + \frac {4(b+\varepsilon)^{2} d_{1} \sigma^{2}}{(C +o(1))\varepsilon^{2}} \right) \ol{G}(x) \dd{P}(X>y).
\end{align*}
Therefore, 
\begin{align*}
  \limsup_{x \to \infty} \frac 1{\ol{G}(x)} \dd{P}\left( A + \sum_{i=1}^{X} B_{i} > x, y < X \le d_{1} x \right) \le \left(c_1\delta(\varepsilon) + \frac{4(b+\varepsilon)^{2} d_{1} \sigma^{2}} {C \varepsilon^{2}} \right) \dd{P}(X>y),
\end{align*}
where $ \dd{P}(X>y)$ in the RHS can be made arbitrarily small by taking 
$y$ sufficiently large. Thus, the numerator in the ratio is ignorable with respect to $\overline{G}(x)$, and we can conclude that the tail asymptotics coincide with the lower bound.

In the case (II), we have
\begin{align}
\label{eq:I+ 2a}
 & \quad \dd{P}\left( A + \sum_{i=1}^{X} B_{i} > x, y < X \le d_{1} x \right) \nonumber\\
 & \quad \le \dd{P}\left(A + \max_{1 \le n \le X} S_{n} > \frac {\varepsilon/2}{b+\varepsilon} x, X > y \right) \nonumber\\
 & \quad \le \dd{P}\left(A + \max_{n \ge 1} S_{n} > \frac {\varepsilon/2}{b+\varepsilon} x \right) \dd{P}(X > y).
\end{align}
Since $H_{I}$ is subexponential, the integrated tail distribution of $B_{i} - (b+\varepsilon/2)$ is also subexponential. Since $S_{n}$ has a negative mean drift, we have, from Theorem 5.2 of \cite{FossKorsZach2011} and \eq{c2 additional},
\begin{align*}
  \dd{P}\left(\max_{n \ge 1} S_{n} > x \right) \sim \frac 1{b+\varepsilon/2} \ol{H}_{I}(x) \lesssim \frac 1{b+\varepsilon/2} \ol{G}(x).
\end{align*}
Hence, similar to \eq{I+ 2aa}, it follows from \eq{I+ 2a} that
\begin{align}
\label{eq:I+ 2b}
 & \dd{P}\left( A + \sum_{i=1}^{X} B_{i} > x, y < X \le d_{1} x \right) \lesssim \left(c_{1} + \frac 1{b+\varepsilon/2}\right) \ol{G}\left(\frac {\varepsilon/4}{b+\varepsilon} x \right) \dd{P}(X > y).
\end{align}
Thus, this term is ignorable compared with $\ol{G}(x)$ by $G \in \sr{D}$ since $\dd{P}(X > y)$ can be arbitrarily small for large $y$, while $\ol{G}\left(\frac {\varepsilon/2}{b+\varepsilon} x \right)/\ol{G}(x)$ is bounded. This completes the proof for the case (II).

\section{Two extensions of the model}
\label{sect:extended}
\setnewcounter

\subsection{Continuous state version}
\label{sect:continuous}

A natural continuous counterpart of \eq{fixed point 1} is
\begin{align}
\label{eq:fixed point 3}
  \widetilde{X} =_{st} \widetilde{A} + \int_{0}^{\widetilde{X}} d\widetilde{B}(t),
\end{align}
where $\widetilde{A}$ is a non-negative random variable independent of $\widetilde{X}$, and $\widetilde{B}(t)$ is a non-decreasing process with stationary independent increments which is independent of $\widetilde{A}$ and $\widetilde{X}$. That is, $\widetilde{B}(\cdot)$ is a non-decreasing Levy process
(subordinator).

We consider a simple case that $\widetilde{B}(\cdot)$ is a compound Poisson process. Namely,
\begin{align}
\label{eq:B tilde 1}
  \widetilde{B}(t) = \sum_{i=1}^{N(t)} \widetilde{B}_{i}, \qquad t \ge 0,
\end{align}
where $N(t)$ is the Poisson process with intensity $\lambda > 0$, and $\widetilde{B}_{i}$ for $i=1,2,\ldots$ are non-negative $i.i.d.$ random variables which are independent of $N(\cdot)$. Assume that $\widetilde{A}$ is a non-negative random variable independent of everything else. Then, \eq{fixed point 3} becomes
\begin{align}
\label{eq:fixed point 4}
  \widetilde{X} =_{st} \widetilde{A} + \sum_{i=1}^{N(\widetilde{X})} \widetilde{B}_{i}.
\end{align}

Similar to the Markov chain $\{X_{n}\}$ defined by \eq{Markov X 1}, we recursively define a discrete time Markov process $\widetilde{X}_{n}$ with state space $\dd{R}_{+}$. This model may be applied, say, to an energy reproduction system. In this system, $\widetilde{A}$ is a base production of energy, and extra energy is reproduced according to the compound Poisson process in the time interval whose length equals the amount of the previous energy production.

The fixed point equation \eq{fixed point 4} can be solved essentially in the same way as \thr{X 1} because
\begin{align*}
  \dd{P}(N(\widetilde{X}) > x) \sim \dd{P}( \lambda\widetilde{X} > x), \qquad x \to \infty.
\end{align*}
Assume that $\widetilde{a} \equiv \dd{E}(\widetilde{A})$ and $\widetilde{b} \equiv \dd{E}(\widetilde{B})$ are finite. Then, \thr{X 1} holds true for $a = \widetilde{a}$ and $b = \lambda \widetilde{b}$ if the solution $\widetilde{X}$ of \eq{fixed point 4} uniquely exists in distribution, where $A$ and $B$ are replaced by $\widetilde{A}$ and $\widetilde{B}$ in the cases (i)--(iii).

\subsection{2nd order branching process with immigration}
\label{sect:2nd order}

In this Subsection, we introduce another extension of the model, formulate a particular result and make short comments on its proof.

Consider a branching process in which two subsequent generations produce the next generation. Namely, let $X_{n}$ be the population of the $n$'th generation, then
\begin{align}
\label{eq:}
  X_{n} = A_{n} + \sum_{i=1}^{X_{n-1}} B_{1,n,i} + \sum_{i=1}^{X_{n-2}} B_{2,n,i}, \qquad n \ge 1.
\end{align}
where $\{A_{n};n \ge 0\}, \{B_{1,n,i}; n \ge 0, i \ge 1\}, \{B_{2,n,i}; n \ge 0, i \ge 1\},$ are sequences of $i.i.d$ non-negative integer-valued random variables, they are mutually independent, and they are independent of $X_{n-1}, X_{n-2}$. We refer to $\{X_{n}; n \ge 0\}$ as a second order branching process.

Let $a = \dd{E}(A)$ and $b_{k} = \dd{E}(B_{k})$ for $k=1,2$. We assume that both $a$ and $b_{k}$ are finite. Then, it is not difficult to see that the process $\{X_{n}; n \ge 0\}$ is stable if and only if $b_{1}+b_{2} < 1$. We assume this stability condition, and consider the following fixed point equation.
\begin{align}
\label{eq:2nd fixed point}
  \left(\begin{array}{c} X \\ Y \end{array}\right) =_{st} \left(\begin{array}{c} A + \sum_{i=1}^{X} B_{1,i} + \sum_{i=1}^{Y} B_{2,i} \\X\end{array}\right).
\end{align}
This fixed point equation uniquely determines the stationary distribution of $X_{n}$ similarly to \eq{fixed point 1}. Note that \eq{2nd fixed point} is equivalent to
\begin{align}
\label{eq:2nd joint distribution}
  \dd{P}(X>x, Y>y) = \dd{P}\left(A + \sum_{i=1}^{X} B_{1,i} + \sum_{i=1}^{Y} B_{2,i} > x, X > y\right).
\end{align}

However, the tail asymptotics for two-dimensional distribution is generally hard to study. So, we restrict our attention to the tail asymptotics
of the linear combination $X+\delta Y$ of $X$ and $Y$, for a particular choice of coefficient $\delta$.  

Note that 
\begin{align}
\label{eq:2nd marginal}
 & \dd{P}(X > x) = \dd{P}(Y > x) = \dd{P}\left(A + \sum_{i=1}^{X} B_{1,i} + \sum_{i=1}^{Y} B_{2,i} > x\right).
\end{align}

From \eq{2nd marginal}, one can find the expectation $m \equiv \dd{E}(X)$:
\begin{align*}
    m = a + (b_{1} + b_{2}) m,
\end{align*}
and therefore, under the stability assumption,
\begin{align}
\label{eq:2nd mean}
  m = \frac {a}{1-(b_{1} + b_{2})} < \infty.
\end{align}

From \eq{2nd fixed point}, we have, for a constant $\delta > 0$, 
\begin{align}
\label{eq:2nd fixed point 2}
  X + \delta Y & =_{st} A + \sum_{i=1}^{X} B_{1,i} + \sum_{i=1}^{Y} B_{2,i} + \delta X \nonumber\\
  & \equiv A + \sum_{i=1}^{X} (\delta + B_{1,i}) + \sum_{i=1}^{Y} B_{2,i}.
\end{align}

In the Proposition below, we provide the distributional tail asymptotics for 
$X+\delta Y$, for a particular choice of $\delta$, under a version of 
condition \eq{G 22}. Under weaker assumptions, one can obtain also upper and
lower bounds.

\begin{proposition}\label{prop2}
%({\bf to be verified!})\\ 
Assume that $a<\infty$ and that $b_1+b_2<1$ (this is the stability condition). 
Assume that there is a reference distribution $G$ such that
\begin{align}\label{eq:c1c3}
  \lim_{x \to \infty} \frac {\dd{P}(A>x)} {\ol{G}(x)} = c_{1}, \qquad \lim_{x \to \infty} \frac {\dd{P}(B_{1,i}>x)} {\ol{G}(x)} = c_{2},
 \qquad \lim_{x \to \infty} \frac {\dd{P}(B_{2,i}>x)} {\ol{G}(x)} = c_{3}   
\end{align}
for some constants $c_{1},c_{2},c_3 \ge 0$ such that $c_{1}+c_{2}+c_3 >0$. 

(I). Let $\delta >0$ be the solution to equation
$
\delta = b_2/(b_1+\delta),
$
i.e. 
$$
\delta = \left(\sqrt{b_1^2 +4b_2}-b_1\right)/2.
$$
Then $b_1+\delta <1$.

(II). Assume that condition \eq{G 22} holds with $c_0=b_1+\delta$. Then 
\begin{align*}
\dd{P} (x<X+\delta Y \le x/(b_1+\delta)) \sim 
(c_1+ \dd{E} X (c_2+c_3))\ol{G}(x)
\end{align*}
and, therefore,
\begin{align*}
\dd{P}(X+\delta Y >x) \sim (c_1+ \dd{E} X (c_2+c_3)) T_{(b_1+\delta)^{-1}}(x).
\end{align*}
%
%\begin{align}
%\label{eq:2nd marginal 1}
%  X + \delta Y & \simeq A + (b_{1} + \delta) X + b_{2} Y + \sum_{i=1}^{\widehat X} %(\delta + B_{i,1}) + \sum_{j=1}^{\widehat Y} B_{j,1} \nonumber\\
%  & \simeq A + (b_{1} + \delta) X + b_{2} Y + \sum_{i=1}^{\widehat{X}} (\delta + %B_{i,1} + B_{i,2}) \nonumber\\
%  & \simeq 
% A + (b_{1} + \delta) \left( X + \frac{b_{2}}{b_1+\delta} Y \right) + \sum_{i=1}%^{\widehat{X}} (B_{i,1} + B_{i,2})\nonumber\\ 
% & =
% (b_{1} + \delta)( X + \delta Y) + A + \sum_{i=1}^{\widehat{X}} (B_{i,1} + B_{i,2})
%  .
%\end{align}
%Here $\widehat{X}=\widehat{Y} = \dd{E} X$ if $\dd{E} X$
%is an integer; otherwise they are two i.i.d. r.v.'s that are independent of
%everything else and taking values $\dd{E} X$ with probability $1-p$ and 
%$\dd{E}X +1$ with probability $p$, where $p= \dd{E}X - [\dd{E}X]$ (therefore,
%$\dd{E} \widehat{X} = \dd{E} X$ again). 
\end{proposition}

Comments on the {\sc Proof} of the proposition. Statement (I) is straightforward.
To obtain the tail asymptotics, we follow the lines of the proof of Theorem 2.1,
with minor modifications. Therefore we replace most of the proof by its sketch,
with providing some details. 

We take the event that the right-hand side of \eq{2nd fixed point 2} exceeds
level $x$, and consider the probabilities $I^-(x,y)$ and $I^+(x,y)$ of the intersection of this event with events $\{X+\delta Y \le y\}$ and  $\{X+\delta Y > y\}$, respectively. For the probability of the first intersection of events, we use again
the result from Proposition \ref{prop1}, while for the second probability we consider again the upper and lower bounds. There are slightly novel arguments
in getting the upper bound only, so we give it in full. We take $\varepsilon \in (0, 1-b_1-\delta)$ and let $d_1=(b_1+\delta+\varepsilon)^{-1}$, $\varepsilon_1 = \varepsilon/2$ and 
$\varepsilon_2=\varepsilon_1b_2/(b_1+\delta)$. We 
have
\begin{align*}
I^+(x,y) 
& =
\dd{P} (A+\sum_1^X (\delta +B_{1,i})+\sum_1^Y B_{2,i} >x, \ X+\delta Y >y)\\
& \le
\dd{P} (X+\delta Y>d_1x) +
 \dd{P} (A+\sum_1^X (\delta +B_{1,i})+\sum_1^Y B_{2,i} >x, \ d_1\ge X+\delta Y >y)\\
 & \equiv
 \dd{P} (X+\delta Y>d_1x) + P(x,y).
 \end{align*}
 Here 
 \begin{align*}
 P(x,y)
 & =
 \dd{P}(A+\sum_1^X (B_{1,i}-b_1- \varepsilon_1) + \sum_1^Y (B_{2,i}-b_2-\varepsilon_2) + (\delta+b_1+\varepsilon_1)X + (b_2+\varepsilon_2)Y>x,\\
 & \quad \quad y<X+\delta Y\le d_1x)\\
 &\le 
 \dd{P}(A+ M_{1,X} + M_{2,Y} + (\delta + b_1+\varepsilon_1)d_1x >x, \ 
 \max (X,Y) > y/2\delta)
 \\
 &\le 
 \dd{P}(A+ M_{1,Z} + M_{2,Z} >\frac{\varepsilon x/2}{\delta + b_1+\varepsilon}, \ 
 Z > y/2\delta)
 \end{align*}
 where, for $n=1,2,\ldots$, 
 $M_{1,n} = \max_{1\le j \le n} \sum_1^j (B_{1,i}-b_1- \varepsilon_1)$,
 $M_{2,n}= \max_{1\le j \le n} \sum_1^j (B_{2,i}-b_2-\varepsilon_2)$, and
 $Z=\max (X,Y)$. Let $\gamma = \varepsilon/6(\delta + b_1+\varepsilon)$.
 Then
 the latter probability is not smaller than
 \begin{align*}
 \dd{P}(A >\gamma x) \dd{P}(Z>y/2\delta) + \sum_{k=1}^2
 \dd{P}(Z>y/2\delta) \dd{P} (M_{k,Z}> \gamma x \ | \ Z>y/2\delta)%
 \end{align*}
 Like in the derivation of the upper bound for $I^+(x,y)$ in the proof of Theorem 2.1, we may 
 use Proposition \ref{prop1} and the property of the class $\cal{D}$ to find that each of the conditional probabilities is proportional
 to $\overline{G}(x)$. Therefore the upper bound to $I^+(x,y)$ is of order
 $c\overline{G}(x)$ where coefficient $c$ may be made as small as one wishes, by taking
 $y$ sufficiently large.

\section{Open problem}
\label{sect:OP}

We do not know, whether \eq{X 2} is the only possible asymptotics for $\dd{P}(X>x)$ in the class  $\sr{L}\cap\sr{D}$ of heavy-tailed distributions. To formulate a more precise open problem, we look closer at equation \eq{fixed point 111}. One may, in turn, represent $X^{(2)}$ as $X^{(2)}=A_2+\sum_1^{X^{(3)}} B_{i,2}$, then use the same representation for $X^{(3)}$, etc. As a result, one can obtain the following a.s. representation for $X=X^{(1)}$:
\begin{align*}
X= A_1+ \sum_{i=1}^{A_2}D_{1,i} + \sum_{i=1}^{A_3}D_{2,i}+ \ldots + 
\sum_{i=1}^{A_n}D_{n-1, i}+\ldots
\end{align*}
where, by convention, $\sum_1^0=0$, all random variables on the right are mutually independent, $A_1=A$ and $\{A_i\}$ are i.i.d., $D_{1,1}=_{st} B$,
$D_{2,1} =_{st} \sum_{i=1}^B B_i$ and, for $n=2,3,\ldots$, $D_{n+1,1} =_{st}
\sum_{i=1}^B D_{n,i}$ where all random variables in the right-hand side of each formula are mutually independent. 

Consider a particular "boundary" example, with $\dd{P}(A>x) = (1+x)^{-1}$, for $\gamma \in (0,1)$, and $\dd{P}(B>x) = L(x) (1+x)^{-1}$ where  $L(x) \sim \left( \log x\right)^{-1-\varepsilon}$, $\varepsilon >0$.  Then $a$ is infinite and $\dd{E}\log \max (1,A)$ is finite. Further, $b = \dd{E} B$ is finite and we can make it smaller than 1.

Then one can use Theorem 7 from \cite{DeFoKo2010} to obtain
$\dd{P}(D_{2,1}>x) \sim 2b \dd{P}(B>x)$ and, using the induction argument,
$\dd{P}(D_{n,1}>x) \sim n b^{n-1} \dd{P}(B>x)$, for any $n=2,3,\ldots$. 
Further, using the uniform convergence result in Theorem 2 of \cite{FoPaZa2005},
one can get the asymptotics
\begin{align}\label{eq:boundary}
\dd{P}\left( \sum_{i=1}^{A_{n+1}} D_{n,i}>x\right) \sim \dd{E}
\left(A{\bf 1}(A\le xb^{-n})\right) \cdot nb^{n-1} \dd{P}(B>x) + \dd{P} (A>xb^{-n}).
\end{align} 
We can expect the PSBJ to hold again and formulate the following {\bf conjecture}:
in the example above,
\begin{align*}
\dd{P}(X>x) \sim \dd{P}(A>x) +\sum_{n=1}^{\infty} \dd{P}\left( \sum_1^{A_{n+1}} D_{n,i}>x\right),
\end{align*}
where the asymptotics for each term in the latter sum are given by \eq{boundary}. However, we do not know how to substantiate these
asymptotics.

\appendix

\section*{Appendix}

\renewcommand{\thesubsection}{\Alph{subsection}}

\renewcommand{\thesubsection}{\Alph{subsection}}
\renewcommand{\theequation}{\thesubsection.\arabic{equation}}
\setcounter{equation}{0}

\subsection{Proof of \lem{G22 sufficient}}
\label{app:G22 sufficient}

By \cite[Proposition 2.2.3]{BingGoldTeug1987}, if $c_{+}(\ol{G}) < 0$, for any $d > c_{+}(\ol{G})$ and any $\gamma > 1$, there exists $x_{0}(\gamma,d)$ such that, for each real $c > 1$ and integer $n \ge 0$,
\begin{align}
\label{eq:G bound 1}
  \ol{G}(c^{n} x)/\ol{G}(x) \le \gamma c^{dn}, \qquad \forall x \ge x_{0}(\gamma,d).
\end{align}
We choose $d < 0$. Since $0 < b < 1$, we can choose $c = c_{+1}(\delta)$ for $\delta > 0$ such that $1 < c_{+1}(\delta) < 1/b$ and $\lim_{\delta \downarrow 0} c_{+1}(\delta) = 1/b$. Then, for any $\varepsilon > 0$, there exists $N(\varepsilon)$ such that
\begin{align*}
  \sum_{n=N(\varepsilon)+1}^{\infty} \ol{G}(c_{+1}^{n}(\delta) x)/\ol{G}(x) \le \gamma \sum_{n=N(\varepsilon)+1}^{\infty} c_{+1}^{-dn}(\delta) < \varepsilon, \qquad \forall x \ge x_{0}(A,-d).
\end{align*}
Hence, we have
\begin{align}
\label{eq:G bound 2}
  1 \le \frac {T_{c_{+1}(\delta)}(x)} {T_{1/b}(x)} & \le \frac {\sum_{n=0}^{N(\varepsilon)} \ol{G}(c_{+1}^{n}(\delta) x)/\ol{G}(x) + \varepsilon} {\sum_{n=0}^{N(\varepsilon)} \ol{G}(b^{-n}x)/\ol{G}(x)} \nonumber\\
  & = 1 + \frac {\sum_{n=0}^{N(\varepsilon)} \frac {\ol{G}(b^{-n}x)} {\ol{G}(x)} \left(\frac {\ol{G}(c_{+1}^{n}(\delta) x)} {\ol{G}(b^{-n}x)} - 1\right) + \varepsilon} {\sum_{n=0}^{N(\varepsilon)} \frac {\ol{G}(b^{-n}x)} {\ol{G}(x)}}.
\end{align}
Since $1 \le \sum_{n=0}^{N(\varepsilon)} \frac {\ol{G}(b^{-n}x)} {\ol{G}(x)} \le N(\varepsilon)$ for all $x \ge 0$ and, for each fixed $n \ge 0$,
\begin{align*}
  \lim_{\delta \downarrow 0} \lim_{x \to \infty} \frac {\ol{G}(c_{+1}^{n}(\delta) x)} {\ol{G}(b^{-n}x)} = \lim_{\delta \downarrow 0} \lim_{x \to \infty} \frac {\ol{G}(x)} {\ol{G}((c_{+1}(\delta) b)^{-n} x)} = 1
\end{align*}
by $G \in \sr{IRV}$, taking the limit of \eq{G bound 2} as $x \to \infty$ then as $\delta \downarrow 0$ yield
\begin{align*}
  1 \le \lim_{\delta \downarrow 0} \limsup_{x \to \infty} \frac {T_{c_{+1}(\delta)}(x)} {T_{1/b}(x)} \le 1 + \varepsilon.
\end{align*}
Letting $\varepsilon \downarrow 0$, we have the first equality of \eq{G 22}. The second equality is similarly obtained by choosing $c = c_{+2}(\delta)$ such that $1 < 1/b < c_{+2}(\delta)$ and $\lim_{\delta \downarrow 0} c_{+2}(\delta) = 1/b$.

The remaining parts of this lemma are obvious, so are omitted.

\subsection{Proof of \lem{monotone 1}}
\label{app:G23}

Let $A$ and $\{B_i\}$ be independent copies of $A_n$ and $\{B_{i,n}\}$, that also does not depend on all $\{X_n\}_{n\ge 0}$ which are obtained by \eq{Markov X 1}.
Since $X_0=0$ and $X_1=A_0 \ge 0$ a.s., we clearly get
$X_{0} \le X_{1}$. Then 
\begin{align*}
  X_{1} =_{st} A+ \sum_{i=1}^{X_{0}} B_{i} \le A+ \sum_{i=1}^{X_{1}} B_{i} =_{st} X_{2}.
\end{align*}
Thus, we have $X_{1} \le_{st} X_{2}$. We then can choose $\widetilde{X}_{1}, \widetilde{X}_{2}$ such that they are independent of $A$ and $\{B_i\}$, $\widetilde{X}_{1}, \le \widetilde{X}_{2}$ a.s. and $\widetilde{X}_{\ell} =_{st} X_{\ell}$ for $\ell=1,2$. Hence,
\begin{align*}
  X_{2} =_{st} A+ \sum_{i=1}^{\widetilde{X}_{1}} B_{i} \le A+ \sum_{i=1}^{\widetilde{X}_{2}} B_{i} =_{st} A+ \sum_{i=1}^{X_{2}} B_{i} =_{st} X_{3}.
\end{align*}
One can repeat this induction argument to conclude that $X_{n} \le_{st} X_{n+1}$ for $n \ge 3$. Since \eq{fixed point 1} has the solution $X$ which is unique in distribution, we have
\begin{align*}
  X_{0} \le_{st} X =_{st} A +\sum_{i=1}^{X} B_i
\end{align*}
Then, we can use the same induction argument as above to get $X_{n} \le_{st} X$. Hence, the distribution of $X_{n}$ weakly converges to some proper distribution $\nu$ as $n \to \infty$. Denote a random variable subject to this $\nu$ by $Y$. Since $X_{n+1} =_{st} A +\sum_{i=1}^{X_n} B_i$ implies that $Y =_{st} A+\sum_{i=1}^Y B_i$, we have $Y =_{st} X$ by the uniqueness of the solution of \eq{fixed point 1} in distribution. Thus, $X_{n}$ converges to $X$ in distribution as $n \to \infty$.

Any stochastically non-decreasing sequence of random variables has a (possibly improper) weak limit, call it $X$. We have
$\dd{E}X_{n+1} = a+ b\dd{E}X_n=\ldots = a(1-b^{n+1})/(1-b)$. 
By the monotone convergence theorem,
$\dd{E} X = \lim_{n\to\infty} \dd{E}X_n = a/(1-b)<\infty$ and, in particular, $X$ is finite a.s. 
This completes the proof. 

\subsection{Example of the tail distribution from the class $\sr{ERV}\setminus\sr{RV}$}
\label{app:ERV}

We provide an example of the tail distribution function $g(x)=\overline{G}(x)$ that
is extended regularly varying, but not regularly varying. 

Let $c>1$ and let $1<a_1 < a_2$.
We assume that function $g(x)$ has the ``cycle'' behaviour and define it by
induction. 
At "time" $t_1=1$, we take $g(t_1)=1$, and the first cycle starts.
Given the $n$'th cycle starts at time $t_n$,
we let $u_n = ct_n$ and define
$
g(t) = (t/t_n)^{-a_1}g(t_n)$ for all $t\in (t_n, u_n]$.
Then let $t_{n+1}=cu_n$ and define
$
g(t) = (t/u_n)^{-a_2} g(u_n)$ for all $t\in (u_n,t_{n+1}]$.

\subsection*{Acknowledgements}
The authors are grateful to Charles Goldie for his invaluable comments on the extended regular variations, $\sr{ERV}$ and for suggesting references \cite{Clin1994,DrasSene1976}.

\def\cprime{$'$} \def\cprime{$'$} \def\cprime{$'$} \def\cprime{$'$}
  \def\cprime{$'$} \def\cprime{$'$} \def\cprime{$'$}


\begin{thebibliography}{8}
\expandafter\ifx\csname natexlab\endcsname\relax\def\natexlab#1{#1}\fi
\expandafter\ifx\csname url\endcsname\relax
  \def\url#1{\texttt{#1}}\fi
\expandafter\ifx\csname urlprefix\endcsname\relax\def\urlprefix{URL }\fi
\providecommand{\eprint}[2][]{\url{#2}}

\bibitem[]{AlsmMein2013}
\textsc{Alsmeyer, G. and Meiners, M.} (2013).
\newblock Fixed points of the smoothing transform: two-sided solutions.
\newblock \textit{Probability Theory and Related Fields}, \textbf{155}, 165--199. 

\bibitem[]{AsmuFoss2018}
\textsc{Asmussen, S. and Foss, S.} (2018).
\newblock Regular variation in a fixed-point problem for single-and multiclass branching processes and queues.
\newblock \textit{Advances in Applied Probability}, \textbf{50(A)}, to appear. 

\bibitem[]{BarcBoszPap2018}
\textsc{Barczy, M., Bosze, Z. and Pap, G.} (2018).
\newblock On tail behaviour of stationary second-order Galton-Watson
processes with immigration.
\newblock \url{https://arxiv.org/abs/1801.07931}

\bibitem[{Basrak-Kulik-Palmowski(2013)}]{BasrKuliPalm2013}
\textsc{Basrak, B., Kulik, R. and Palmowski, Z.} (2013).
\newblock Heavy-tailed branching process with immigration.
\newblock \textit{Stochastic Models}, \textbf{29}, 413--434. 

\bibitem[{Bingham et~al.(1987)Bingham, Goldie and Teugels}]{BingGoldTeug1987}
\textsc{Bingham, N.}, \textsc{Goldie, C.} and \textsc{Teugels, J.} (1987).
\newblock \textit{Regular Variation (Encyclopedia of Mathematics and its
  Applications)}.
\newblock Cambridge University Press.

\bibitem[]{Clin1994}
\textsc{Cline, D.S.} (1994).
\newblock Intermediate regular and $\Pi$ variation.
\newblock \textit{Proceedings of the London Mathematical Society}, \textbf{68}, 594--616. 

\bibitem[{DeFoKo(1997)}]{DeFoKo2010}
\textsc{Denisov, D.; Foss, S.; Korshunov, D.} (2010).
\newblock Asymptotics of randomly stopped sums in the presence of heavy tails.
\newblock \textit{Bernoulli}, \textbf{16}, 971--994. 

\bibitem[{Drasin and Seneta (1976)}]{DrasSene1976}
\textsc{Drasin, D. and Seneta, E.} (1976).
\newblock A generalization of slowly varying functions.
\newblock \textit{Proceedings of the American Mathematical Society}, \textbf{96} 470--472.

\bibitem[{Fayolle et~al.(1999)Fayolle, Iasnogorodski and
  Malyshev}]{FayoIasnMaly1999}
\textsc{Fayolle, G.}, \textsc{Iasnogorodski, R.} and \textsc{Malyshev, V.}
  (1999).
\newblock \textit{Random Walks in the Quarter-Plane: Algebraic Methods,
  Boundary Value Problems and Applications}.
\newblock Springer, New York.

\bibitem[{Fayolle et~al.(1995)Fayolle, Malyshev and
  Menshikov}]{FayoMalyMens1995}
\textsc{Fayolle, G.}, \textsc{Malyshev, V.} and \textsc{Menshikov, M.} (1995).
\newblock \textit{Topics in the constructive theory of countable {Markov}
  chains}.
\newblock Cambridge University Press, Cambridge, UK.

\bibitem[{Foss et~al.(2011)Foss, Korshunov and Zachary}]{FossKorsZach2011}
\textsc{Foss, S.}, \textsc{Korshunov, D.} and \textsc{Zachary, S.} (2011).
\newblock \textit{An Introduction to Heavy-Tailed and Subexponential
  Distributions}.
\newblock Springer Series in Operations Research and Financial Engineering,
  Springer.

\bibitem[{Foss and Miyazawa(2018)}]{FossMiya2018}
\textsc{Foss, S.} and \textsc{Miyazawa, M.} (2018).
\newblock Customer sojourn time in {$GI/G/1$} feedback queue in the presence of
  heavy tails.
\newblock \textit{The Journal of Statistical Physics}.
\newblock To appear, \urlprefix\url{https://arxiv.org/abs/1710.10503}.

\bibitem[{Foss-Palmowski-Zachary(2005)}]{FoPaZa2005}
\textsc{Foss, S.; Palmowski, Z.; Zachary, S.} (2005).
\newblock The probability of exceeding a high boundary on a random time interval for a heavy-tailed random walk.
\newblock \textit{Annals of Applied Probability}, \textit{15}, 1936--1957.  

\bibitem[{Foss and Zachary(2003)}]{FossZach2003}
\textsc{Foss, S.} and \textsc{Zachary, S.} (2003).
\newblock The maximum on a random time interval of a random walk with a
  long-tailed increments and negative draft.
\newblock \textit{Annals of Applied Probability}, \textbf{13} 37--53.

\bibitem[{Foster and Williamson(1971)}]{FostWill1971}
\textsc{Foster, J.} and \textsc{Williamson, J.} (1971).
\newblock Limit theorems for the Galton-Watson process with time-dependent
  immigration.
\newblock \textit{Z. Wahrscheinlichkeitstheorie verw. Gebiete}, \textbf{20}
  227--235.

\bibitem[]{JeleOlve2012}
\textsc{Jelenkovi\'{c}, P.R. and Olvera-Cravioto, M.} (2012).
\newblock Implicit renewal theory and power tails on trees.
\newblock \textit{Advances in Applied Probability}, \textbf{44} 528--561.

\bibitem[{Seneta(1971)}]{Sene1971}
\textsc{Seneta, E.} (1971).
\newblock On invariant measures for simple branching processes.
\newblock \textit{Journal of Applied Probability}, \textbf{8} 43--51.

\end{thebibliography}
\end{document}